\documentclass{amsart}
\usepackage{amsthm}
\usepackage{amssymb}
\usepackage{amsmath}
\usepackage{amscd}

\newtheorem {Theorem}{Theorem}
\numberwithin{Theorem}{section}
\newtheorem {Lemma}[Theorem]{Lemma}

\theoremstyle{definition}
\newtheorem{Definition}[Theorem]{Definition}
\theoremstyle{remark}
\newtheorem{Remark}[Theorem]{Remark}

\newtheorem{Assumption}[Theorem]{Assumption}
\numberwithin{equation}{section}
\expandafter\chardef\csname pre amssym.def at\endcsname=\the\catcode`\@ \catcode`\@=11
\def\undefine#1{\let#1\undefined}
\def\newsymbol#1#2#3#4#5{\let\next@\relax
 \ifnum#2=\@ne\let\next@\msafam@\else
 \ifnum#2=\tw@\let\next@\msbfam@\fi\fi
 \mathchardef#1="#3\next@#4#5}
\def\mathhexbox@#1#2#3{\relax
 \ifmmode\mathpalette{}{\m@th\mathchar"#1#2#3}%
 \else\leavevmode\hbox{$\m@th\mathchar"#1#2#3$}\fi}
\def\hexnumber@#1{\ifcase#1 0\or 1\or 2\or 3\or 4\or 5\or 6\or 7\or 8\or
 9\or A\or B\or C\or D\or E\or F\fi}
\font\teneufm=eufm10 \font\seveneufm=eufm7 \font\fiveeufm=eufm5
\newfam\eufmfam
\textfont\eufmfam=\teneufm \scriptfont\eufmfam=\seveneufm
\scriptscriptfont\eufmfam=\fiveeufm

\catcode`\@=\csname pre amssym.def at\endcsname
\newcounter{remark}
\newcommand{\bg}{\begin{equation}}
\newcommand{\ed}{\end{equation}}
\newcommand{\bga}{\begin{eqnarray}}
\newcommand{\eda}{\end{eqnarray}}

\def\cbdu{\hfill{$\Box$}}

\def  \12  {{\frac{1}{2}}}

\begin{document}

\title[Steady-State Navier-Stokes in $\mathbb{R}^3$]{Existence and Stability of Steady-State Solutions with Finite Energy for the Navier-Stokes equation in the Whole Space}

\author[Clayton Bjorland ]{Clayton Bjorland }

\address{Department of Mathematics, UC Santa Cruz, Santa Cruz, CA 95064,
USA}

\author[Maria E. Schonbek ]{Maria E. Schonbek }

\address{Department of Mathematics, UC Santa Cruz, Santa Cruz, CA 95064,
USA}

\thanks{The work of M. Schonbek was partially supported by NSF Grant DMS-0600692}
\thanks{The work of C. Bjorland was partially supported by NSF grant OISE-0630623 and UCSC Chancellor’s Dissertation-Year Fellowship}

\email{cbjorland@math.ucsc.edu} \email{schonbek@math.ucsc.edu}


 \

\date{\today}



\bigskip

\begin{abstract}
We consider the steady-state Navier-Stokes equation in the whole space $\mathbb{R}^3$ driven by  a forcing function $f$.  The class of source functions $f$ under consideration yield the existence of at least one solution with finite Dirichlet integral ($\|\nabla U\|_2<\infty$).  Under the additional assumptions that $f$ is absent of low modes and the ratio of $f$ to viscosity is sufficiently small in a natural norm we construct solutions which have finite energy (finite $L^2$ norm).  These solutions are unique among all solutions with finite energy and finite Dirichlet integral.  The constructed solutions are also shown to be stable in the following sense: If $U$ is such a solution then any viscous, incompressible flow in the whole space, driven by $f$ and starting with finite energy, will return to $U$.
\end{abstract}

\maketitle
\section{Introduction}
The classical theory of viscous, incompressible fluid flow  is governed by the famous Navier-Stokes equations:
\begin{align}\label{NS:PDE}
u_t+u\cdot\nabla u+\nabla p &= \nu \triangle u +f\\
u(0)=u_0 \ \ \ \nabla \cdot u &=0\notag
\end{align}
A large area of modern research is devoted to deducing qualitative properties of solutions for these equations when they are complemented with initial and boundary conditions and certain restraints are placed on $f$ and $\nu$.  The investigations in this subject are too numerous to attempt to list here so we will limit ourselves  to discussion directly related to the topic of this paper: the steady state Navier-Stokes equation in the whole space $\mathbb{R}^3$.  

A steady state (sometimes called stationary in the literature) solution $U$ of the Navier-Stokes equation is one for which  $\partial_t U=0$, that is the solution is constant in time.  Such solutions solve the following PDE, which will be our main point of investigation.
\begin{align}\label{stationary:PDE}
U\cdot\nabla U +\nabla p &=\nu\triangle U +f\\
\nabla\cdot U &=0\notag
\end{align}
For our purposes this PDE is supplemented with the idea that $U$ tends to zero as $|x|$ becomes large, made precise by working in functions spaces which are completions of smooth functions with compact support.

Roughly speaking, the investigation of steady state solutions can be broken into two regimes: bounded and unbounded domains.  In the former situation much progress has been made using  a Poincar\'e type inequality ($\|U\|_2\leq C\|\nabla U\|_2$) to deduce quickly that solutions have finite energy.  In the case where there is no Poincar\'e inequality it is desirable to to find conditions on $f$ which will guarantee, a priori, finite energy of a solution.  One of the benefits of establishing existence of solutions with finite energy in unbounded domains is that many of the techniques developed using the Poincar\'e inequality can be applied, but the consideration of unbounded domains is not a needless complication.  Indeed, many physical problems are best stated in the whole space or in exterior domains where there is no Poincar\'e inequality.  Moreover, the situation in the whole space is theoretically important as Leray observed in his seminal paper outlining modern analysis of the Navier-Stokes equations \cite{MR1555394}:

\begin{quote}\textit{L'absence de parois introduit certes quelques complications concernant l'allure \`a l'infini des fonctions inconnues, mais simplifie beaucoup l'expos\'e et met mieux en lumigre les difficult\'es essentielles;}
\end{quote}

The main goal of this paper is to develop a new technique which will allow, with certain conditions on $f$, the construction of solutions for the steady state Navier-Stokes equation in the whole space with finite energy.  The assumptions we impose on $f$ limit the amount of \textit{low frequency} information and require that the ratio of $f$ to $\nu$ is small in a natural norm.  Once we have established the finite energy of solutions we deduce uniqueness in the class of solutions with finite energy and prove these solutions are stable in a strong sense referred to in the literature as \textit{nonlinearly stable}.  In other situations (not the whole space) non-uniqueness for solutions of (\ref{stationary:PDE}) has been demonstrated for solutions with $f$ large compared with $\nu$ (see \cite{MR1846644}, Chapter 2, and references therein) so we suspect the smallness assumption we make on $f$ is necessary and natural for this result.  It is currently unknown if the assumption on low frequency information is natural or a byproduct of our technique.

\bigskip


\subsection{Statement of Results}
Modern analysis of the steady state Navier-Stokes equation in unbounded domains can be traced back to \cite{LerayII} and \cite{Oseen}, these ideas were further developed in \cite{MR0107442}, \cite{MR0166498}, and \cite{MR0182816} which work with the notion of a \textit{physically reasonable solution}.  The authors were concerned with the physically interesting problem of solutions in exterior domains of which the whole space is a special case.  We summarize (perhaps too succinctly) these works with the following idea:  If $f$ is such that the Dirichlet integral $\|\nabla U\|_2$ is finite ($|x|f\in L^2$ is sufficient) then there exists a unique \textit{physically reasonable solution} $U$ in an exterior domain.  This solution is physically reasonable in the sense that it approaches a constant (possibly non-zero) vector field like $|x|^{-1}$ as $|x|$ becomes large and the uniqueness is among all such functions.  These ideas were expanded further in \cite{MR1158939}, \cite{MR1351321}, \cite{MR1467557}, \cite{MR2299756}, \cite{MR0412639}, and \cite{MR1360998}.  The methods rely heavily upon analysis of the Green's function for the domain in question and are quite different from the approach presented in this paper.  

Our construction of solutions with finite energy is based on a well known formal observation: if $\Phi$ is the fundamental solution for the heat equation then $\int_0^\infty \Phi(t,\cdot)\,dt$ is the fundamental solution for Poisson's equation.  Using this idea it is possible to make a time dependent PDE similar to the Navier-Stokes equation with $f$ as initial data with a solution that can be formally integrated in time to find a solution of (\ref{stationary:PDE}).  At this point our analysis turns to the theory of energy decay for fluid equations.  If the decay of the new time dependent PDE is fast enough (an integral over all time converges) we can deduce a finite energy bound for (\ref{stationary:PDE}).  It is known that the decay rate of solutions for parabolic PDEs in the whole space is intimately related to the shape of the initial data near the origin in Fourier space; the low frequency assumption we make on $f$ is enough to guarantee the convergence of the required time integrals.   This idea is further outlined in Section 2 and made precise in Section 3.  

A particularly useful technique for estimating energy decay is the \textit{Fourier Splitting Method}  which was used in \cite{MR775190} to establish energy decay for initial data $u_0\in L^2\cap L^1$ and later for initial data $u_0\in L^2$ in \cite{MR1432588}.  Other works in this area include \cite{MR1158939},  \cite{MR835398}, \cite{MR760047}, \cite{MR1360998}, \cite{MR1815476}, \cite{MR968313}, \cite{MR837929}, \cite{MR1312701},  \cite{MR1776344}, \cite{MR1396285}, \cite{MR881519}, and \cite{MR1896923}.  In essence we are trading bounds on the Green's function for energy decay theorems which we base on the Fourier Splitting Method.  The assumption on low frequency information is stronger then the classical assumptions but our conclusion is stronger.  Of course the previous results consider the more complicated cases of external domains which are not handled within but we hope that with decay theorems for external domains one can use the technique presented here to obtain similar results.  We now state precisely the main theorem proved.  In the following statement $\mathring{H}^1_\sigma$ is the completion of smooth divergence free functions of compact support under the norm $\|\nabla \cdot\|_2$ and $H^1_\sigma=\mathring{H}^1_\sigma\cap L^2_\sigma$.  Also, $X=(\mathring{H}_\sigma^1)'\cap L^2_\sigma$.  The requirement $f\in X$ implies the classical assumptions $f\in L^2$ and $\|\nabla U\|_2 <\infty$.  The later is known as a \textit{finite Dirichlet integral} and is sometimes implied by the restriction $|x|f\in L^2$ in the literature.

\begin{Theorem}\label{L2norm:theorem}
Let $M>0$ and $f\in X$ satisfy the following assumption: 
\begin{itemize}
\item[(A)]There exists a $\rho_0$ such that $\hat{f}(\xi)=0$ for almost every $|\xi|<\rho_0$
\end{itemize}  
Then there exists a constant $C(\rho_0,\nu,M)$ so that if $\|f\|_X\leq C(\rho_0,\nu,M)$ the following hold:
\begin{itemize}
\item[(i)] The PDE (\ref{stationary:PDE}) has a weak solution $U\in H^1_\sigma$.  It is a weak solution in the sense that for any divergence free function of compact support $\phi$, 
\begin{equation}\label{weakstatement:stationary} 
 <U\cdot\nabla U,\phi> + \nu<\nabla U, \nabla \phi> = <f,\phi>
\end{equation}
\item[(ii)] This solution satisfies $\|U\|_2\leq M$ and $\|\nabla U\|_2\leq \nu^{-1}\|f\|_X$.
\item[(iii)] This solution is unique among all solutions which have a finite $L^2$ norm and satisfy $\|\nabla U\|_2\leq \nu^{-1} \|f\|_X$.
\end{itemize}
\end{Theorem}
\begin{Remark}
The behavior of the constant $C(\rho_0,\nu,M)$ allows large $f$ when the Reynolds number is small (see remark \ref{ul2M:remark}).
In this work we assume the Fourier transform of $f$ is zero in some neighborhood of the origin, this corresponds to exponential decay for the heat flow starting with initial data $f$.  It is possible to relax the hypothesis so the heat flow is algebraic but fast and not significantly change any of the proofs presented here.  This could be accomplished by requiring $f$ to \textit{behave like} the polynomial $|\xi|^p$ near the origin in Fourier space where $p$ is some sufficiently large number (see \cite{BSII}).
\end{Remark}

When a solution is known to have finite energy the well developed energy stability arguments can be applied when $f$ is small measured against $\nu$, see \cite{MR1325465}, \cite{MR2098531}, \cite{MR0105250},  and \cite{MR1140924}.  This is argued by showing solutions of a nonlinear parabolic PDE, similar to the Navier-Stokes equation and found by subtracting the steady state, tends zero as time becomes large.  Again we turn to energy decay methods, specifically a method applied to the Navier-Stokes equation in \cite{MR1432588}, to show the decay.  The method consists of estimating the high and low frequencies of the solution separately and the estimates of the high frequency rely on the Fourier Splitting Method.  Section 4 is dedicated to establishing the following theorem.

\begin{Theorem}
Let $f$ satisfy the assumptions of Theorem \ref{L2norm:theorem} and be such that $\|f\|_X$ is less then the constant given by the theorem.  There exists another constant $C(\nu)$ such that $\|f\|_X\leq C(\nu)$ implies $U$ is stable in the following sense: if $w_0\in L^2_\sigma$ is a perturbation and $u$ is a solution of the Navier-Stokes equation (\ref{NS:PDE}) with initial data $w_0+U$ which satisfies, for any $T>0$,
\begin{equation}
u\in L^{\infty}(0,T;L^2)\cap L^{2}(0,T;\mathring{H}^1_\sigma) \notag
\end{equation}
then
\begin{itemize}
 \item[(i)] For every $\epsilon >0$ there is a $\delta>0$ such that 
\begin{equation}\notag
 \|w_0\|_2\leq \delta\ \ \  \mathrm{implies} \ \ \ \sup_{t\in\mathbb{R}^+}\|u(t)- U\|_2 \leq \epsilon
\end{equation}
 \item[(ii)]$u(t)$ tends to $U$ as time becomes large, that is
\begin{equation}\notag
 \lim_{t\rightarrow\infty}\|u(t)-U\|_2=0
\end{equation}
\end{itemize}
\end{Theorem}

\bigskip

\subsection{Notation}
Unless otherwise noted all integrals in this paper are taken over the whole space $\mathbb{R}^3$, $C_0^\infty$ denotes the space of smooth functions with compact support.  
\begin{align}
&<f,g> = \int f\cdot g&
&\mathcal{V}=\{ \phi \in C_0^\infty | \nabla\cdot \phi=0\}\notag\\
&\|\cdot\|_p=\left(\int |\cdot|^p \right)^{1/p}&
&L^p_\sigma =\{\textrm{ completion of $\mathcal{V}$ under the norm $\|\cdot\|_p$}\}\notag\\
&\|\cdot\|_{\mathring{H}^1}=\|\nabla \cdot\|_2&
&\mathring{H}^1_\sigma = \{\textrm{ completion of $\mathcal{V}$ under the norm $\|\cdot\|_{\mathring{H}^1}$}\}\notag\\
&H^1_\sigma = L^2_\sigma\cap \mathring{H}^1_\sigma&
&(\mathring{H}^1_\sigma)' = \{\textrm{dual of $\mathring{H}^1_\sigma$}\}\notag\\
&X=L^2_\sigma\cap (\mathring{H}^1_\sigma)'&
&\|\cdot\|_X=\max\{\|\cdot\|_2,\|\cdot\|_{(\mathring{H}^1_\sigma)'}\}\notag\\
&\hat{f}(\xi) = \int f(x)e^{-2\pi i x\cdot\xi}\, dx&
&\check{f}(x) = \int f(\xi)e^{2\pi i x\cdot\xi}\, d\xi\notag
\end{align}
We will typically write an element $f\in (\mathring{H}^1_\sigma)'$ as ``$f$'' when we really mean the map ``$\phi \rightarrow <f,\phi>$''.  To denote general constants we use $C$ which may change from line to line.  In certain cases we will write $C(\alpha)$ to emphasize the constants dependence on $\alpha$.  In a similar way we write $p$ to denote general potentials (used to describe the pressure, one instance of $p$ may not be the same as another even on neighboring lines).  The variable $\xi$ is reserved for working in Fourier Space.  

\bigskip

\section{Preliminaries}
Existence of weak solutions for (\ref{stationary:PDE}) is well known, see for example \cite{MR972259}, \cite{MR0107442}, \cite{MR0132307}, \cite{MR0050423},  \cite{MR0254401}, \cite{LerayII}, \cite{MR1846644}.  A typical approach to constructing weak solutions for this PDE is to construct approximations with the Galerkin Method and use \textit{a priori} bounds with the Banach-Alaoglu theorem to find a subsequence of approximations converging weakly to a possible solution.  Some stronger compactness property is then used to pass the sequence through the nonlinear term and establish the limit is indeed a solution.  A good \text{a priori} bound for this approach, and a bound we will rely on throughout is:
\begin{equation}\label{stationary:bound}
\|\nabla U\|_2^2\leq \nu^{-2}\|f\|_X^2
\end{equation}
This is essentially the classical assumption that $U$ has a finite Dirichlet integral but we derive it from our assumption $f\in X$ using the estimate
\begin{align}
|<f,U>|\leq \|f\|_X\|\nabla U\|_2\notag
\end{align}
The bound (\ref{stationary:bound}) is proved formally by multiplying (\ref{stationary:PDE}) by $U$ then integrating by parts.  Noting the specific form of the nonlinearity,
\begin{equation}\label{bilinear:relation}
<w\cdot\nabla U, U>=0
\end{equation}
This relation holds when $\nabla\cdot w=0$ and the integral is absolutely summable, it can be proved for functions of compact support using integration by parts then extended to other classes of functions with a density argument.  It holds in three dimensions when $w\in L^3_\sigma$ and $U\in \mathring{H}^1_\sigma$, or when $U\in L^3_\sigma\cap \mathring{H}^1_\sigma$ and $w\in \mathring{H}^1_\sigma$ since either assumption implies summability.  

Fix $f$ and $U$ as a solution to (\ref{stationary:PDE}) ($U$ does \textit{not} depend on time), we would like to find conditions on $f$ which guarantee $\|U\|_2<\infty$.  One of the key steps of our approach is to establish ``fast'' decay of solutions to the system
\begin{align}
 v_t +U\cdot\nabla v +\nabla p &= \nu\triangle v\label{stationaryneighbor:PDE}\\
v(0)=f \ \ \ \nabla\cdot v&=0\notag
\end{align}
Formally, if $v$ is a solution of (\ref{stationaryneighbor:PDE}) then $\tilde{U}=\int_0^\infty v(t)\, dt$ solves
\begin{align}
 U\cdot\nabla \tilde{U} +\nabla p &= \nu\triangle \tilde{U} +f\notag\\
 \nabla\cdot \tilde{U}&=0\notag
\end{align}
Recall we have fixed $U$ earlier and it is also a solution for this PDE since it satisfies (\ref{stationary:PDE}).  As this PDE is linear and $\nabla\cdot U=0$, solutions are unique and we may conclude $\tilde{U}=U$.  Using Minkowski's Inequality for integrals we can see how the $L^2$ decay of $v$ relates to the $L^2$ norm of $U$:
\begin{equation}
 \|U\|_2 = \|\int_0^\infty v(t) dt\|_2 \leq \int_0^\infty \|v(t)\|_2 dt\notag
\end{equation}
In summary, if $\|v(t)\|_2\leq C(1+t)^{-\beta}$ with $\beta >1$ we can expect $U\in L^2$.

Through a standard Fourier splitting argument we can only hope 
\begin{align}
\|v\|_2\leq C(1+t)^{-3/4}\notag
\end{align} 
where the hold up for faster decay is the initial data.  To get around this problem we will measure the difference 
\begin{align}
w=v-\Phi\ \ \ \  \mathrm{where} \ \ \ \ \Phi=e^{\nu\triangle t}f\notag
\end{align}
Here $\Phi$ is the solution to the heat equation with initial data $f$.  The function $w$ satisfies a parabolic equation with zero initial data and a forcing term which we can control by restricting $f$:
\begin{align}\label{difference:PDE}
 w_t +U\cdot\nabla w +\nabla p &= \nu\triangle w - U\cdot\nabla \Phi\\
\nabla\cdot w =0 \ \ \ w(0) &= 0\notag
\end{align}
We can expect $\|w(t)\|_2$ to decay as $(1+t)^{-5/4}$ and if the heat flow corresponding to $f$ decays at least as fast we can say the same about $v$, to make other parts of the argument work we need $\Phi$ to decay faster.  It is well known that the energy decay of the heat flow corresponding to $f$ is intimately related to the behavior of $\hat{f}$ near the origin, therefore an assumption made on the decay of $\Phi$ is really an assumption on $\hat{f}$ near the origin.  
With this in mind make the following assumption on $f$:
\begin{Assumption}\label{assumption1:ass}
$f\in X$ and there exists a $\rho_0>0$ such that $\hat{f}(\xi)=0$ for every $|\xi|<\rho_0$.
\end{Assumption}
\begin{Remark}
This assumption is really a bandpass filter for $f$ which eliminates low frequencies.  It is known that the corresponding heat energy decays exponentially, a fact demonstrated in the following lemma.  Strictly speaking, one can relax the assumption on $f$ so the heat energy decays at an algebraic (not exponential) rate and use the same method outlined in this paper, such $f$ will need to be behave like a polynomial $|\xi|^p$ near the origin in Fourier space.  See \cite{BSII}.
\end{Remark}

\begin{Lemma}\label{assumption1:lemma}
If $f$ satisfies Assumption \ref{assumption1:ass} and $\Phi=e^{\nu\triangle t}f$, then
\begin{equation}\label{assumption1:bound}
\|\Phi\|_2^2\leq e^{-2\nu\rho_0t}\|\hat{f}\|_2^2
\end{equation}
\end{Lemma}
\begin{proof}
The proof is quickly checked using the bound 
\begin{align}
|\hat{\Phi}|&=|e^{-\nu|\xi|^2t}\hat{f}|\notag\\
&\leq  e^{-\nu\rho_0t}|\hat{f}|\notag
\end{align}
and computing the $L^2$ norm with the aid of the Plancherel theorem.
\end{proof}

\bigskip

\section{$L^2$ Bounds for Stationary Solutions of the NSE}
Throughout this section we will assume $f$ satisfies Assumption \ref{assumption1:ass} and therefore $\Phi=e^{\nu\triangle t}f$ satisfies (\ref{assumption1:bound}).  We are focused on the study of solutions for the two auxiliary PDEs:
\begin{align}\label{stationary:aPDE}
U^i\cdot\nabla U^{i+1} +\nabla p &=\nu\triangle U^{i+1} +f\\
\nabla\cdot U^{i+1} &=0\notag
\end{align}
and
\begin{align}\label{difference:aPDE}
 w^{i+1}_t +U^{i}\cdot\nabla w^{i+1} +\nabla p &= \nu\triangle w^{i+1} - U^{i}\cdot\nabla \Phi\\
\nabla\cdot w^{i+1} =0 \ \ \ w^{i+1}(0) &= 0\notag
\end{align}
When dealing with either PDE we take the function $U^i\in H^1_\sigma$ fixed before hand.  These PDE's will be used recursively to find approximate solutions for (\ref{stationary:PDE}) and (\ref{difference:PDE}) respectively.  In subsection 3.1 we recall existence theorems for these equations and Subsection 3.2 contains the decay rate calculations for $w^{i+1}$.  In Subsection 3.3 we make precise the notion $U^{i}=\int_0^{\infty}v^{i}(t)\,dt$ which is then combined with decay calculations in Subsection 3.4 to find uniform bounds on $U^{i}$ and show it is a Cauchy sequence in $\mathring{H}^1_\sigma$ whose limit is a solution of (\ref{stationary:PDE}).

\bigskip

\subsection{Existence Theorems}

\begin{Theorem}\label{stationary:aexist}
Let $U^i\in H^1_\sigma$ and $f\in X$.  There exists a unique weak solution $U^{i+1}$ to the PDE (\ref{stationary:aPDE}) in the sense that for any $\phi\in \mathcal{V}$,
\begin{align}
 <U^i\cdot\nabla U^{i+1},\phi> + \nu<\nabla U^{i+1}, \nabla \phi> = <f,\phi>
\end{align}
Moreover, this solution satisfies 
\begin{equation}\label{stationary:abound}
\|\nabla U^{i+1}\|_2^2\leq \nu^{-2}\|f\|_X^2 
\end{equation}
\end{Theorem}
\begin{proof}
We only outline the proof as similar PDEs are solved in the literature, see \cite{MR972259}, \cite{MR0107442}, \cite{MR0132307}, \cite{MR0050423},  \cite{MR0254401}, \cite{LerayII}, and \cite{MR1846644}.  A typical approach is to construct Galerkin approximations $\{U_n^{i+1}\}_{n\in\mathbb{N}}$ by projecting the PDE onto finite dimensional subspaces of $H^1_\sigma$.  A uniform bound similar to (\ref{stationary:abound}) can be proved for each Galerkin approximation using an argument similar to that following (\ref{stationary:bound}).  Once this bound is established it is possible to use the Banach-Alaoglu Theorem to find a subsequence $\{U_n^{i+1}\}_{n\in N\subset \mathbb{N}}$ that converges weakly in $\mathring{H}^1_\sigma$.  The weak convergence is enough to pass to a limit in the linear terms.  To pass through the nonlinear term one uses a stronger compactness theorem in the support of the test function $\phi$.
\end{proof}

\begin{Theorem}\label{difference:aexist}
Let $U^i\in H^1_\sigma$ satisfy 
\begin{align}\label{difference:aass}
\|\nabla U^i\|_2\leq \nu^{-1}\|f\|_X
\end{align}
and $f$ satisfy Assumption \ref{assumption1:ass} with $\Phi=e^{\nu\triangle t}f$.  There exists a unique weak solution $w^{i+1}\in L^\infty(\mathbb{R}^+,L^2_\sigma)\cap L^2(\mathbb{R}^+,\mathring{H}^1_\sigma)$ to the PDE (\ref{difference:aPDE}) in the sense that for any $\phi\in C^1(\mathbb{R}^+;\mathcal{V})$,
\begin{align}
<w^{i+1}_t,\phi>+<U^{i}\cdot\nabla w^{i+1},\phi>&=-\nu<\nabla w^{i+1},\nabla\phi>-<U^{i}\cdot\nabla \Phi,\phi>\label{difference:aweakformulation}\\
\nabla\cdot w^{i+1}=0 \ \ \ w^{i+1}(0)&=0\notag
\end{align}
Moreover, this solution satisfies 
\begin{equation}\label{wenergy:apbound}
\sup_t \|w^{i+1}(t)\|_2^2 +\nu\int_0^\infty \|\nabla w^{i+1}(s)\|_2^2\,ds \leq C\rho_0^{-\frac{1}{2}}\nu^{-4}\|f\|_X^4
\end{equation}
\end{Theorem}
\begin{proof}
The PDE in question is closely related to the Navier-Stokes equation and we refer to the literature for similar arguments, see \cite{MR673830}, \cite{MR972259}, \cite{MR589434}, \cite{MR0050423}, \cite{MR1555394}, and \cite{MR1846644}.  It is typical to construct a sequence of Galerkin approximations which satisfies a uniform estimate similar to (\ref{wenergy:apbound}) then use compactness arguments to pass through the limit.  We give now a formal proof of (\ref{wenergy:apbound}) which can be used as an \textit{a priori} estimate in this approach.

Multiply (\ref{difference:PDE}) by $w^{i+1}$ and integrate by parts, then use the bilinear relation (\ref{bilinear:relation}) to find
\begin{align}
\frac{1}{2}\frac{d}{dt}\|w^{i+1}\|_2^2+\nu\|\nabla w^{i+1}\|_2^2 &= <U^{i}\cdot\nabla w^{i+1},\Phi>\notag\\
&= \|U^{i}\|_6\|\nabla w^{i+1}\|_2\|\Phi\|_3\notag\\
&\leq \frac{C}{\nu}\|U^{i}\|^2_6\|\Phi\|^2_3 +\frac{\nu}{2}\|\nabla w^{i+1}\|_2^2\notag
\end{align}
The last line was obtained using H\"older's inequality then Cauchy's inequality.  Putting this together with the Gagliardo-Nirenberg-Sobolev inequality and the assumed bound on $\|\nabla U^{i}\|_2$ yields
\begin{align}
\frac{d}{dt}\|w^{i+1}\|_2^2+\nu\|\nabla w^{i+1}\|_2^2 \leq \frac{C}{\nu^3}\|f\|_X^2\|\Phi\|_3^2\label{temp1}
\end{align}
Using Lemma \ref{assumption1:lemma}, the bound  (\ref{assumption1:bound}) implies
\begin{align}
\int_0^t\|\Phi(s)\|_2^2\,ds \leq \|f\|_X^2\int_0^te^{-2\nu\rho_0s}\,ds\leq \frac{\|f\|_X^2}{2\nu\rho_0}\notag 
\end{align}
Together with the Gagliardo-Nirenberg-Sobolev inequality and the heat property $2\nu\int_0^\infty \|\nabla \Phi(s)\|_2^2\,ds\leq \|f\|^2_2$ we estimate
\begin{align}\label{assumption1:heatl3}
 \int_0^t \|\Phi(s)\|_3^2\, ds &\leq \int_0^t \|\Phi(s)\|_2\|\nabla\Phi(s)\|_2\, ds\notag\\
&\leq \left(\int_0^t \|\Phi(s)\|_2^2\, ds\right)^{\frac{1}{2}}\left(\int_0^t \|\nabla\Phi(s)\|_2^2\, ds\right)^{\frac{1}{2}}\notag\\
&\leq C\rho_0^{-\frac{1}{2}}\nu^{-1}\|f\|_X^2
\end{align}
Integrating (\ref{temp1}) in time then applying (\ref{assumption1:heatl3}) finishes the proof.
\end{proof}

\begin{Remark}\label{justifymultiply:remark}
 In the theorems above the assumption $U^i\in H^1_\sigma$ is enough to ensure $U^{i}\cdot\nabla U^{i+1}\in (\mathring{H}^1_\sigma)'$ and $U^{i}\cdot\nabla w^{i+1}\in (\mathring{H}^1_\sigma)' \ a.e.$  That is 
\begin{equation}\notag
|<U^i\cdot\nabla U^{i+1}, \phi>| \leq C\|\nabla \phi\|_2^2 \ \ \ \mathrm{and} \ \ \ |<U^i\cdot\nabla w^{i+1}, \phi>| \leq C\|\nabla \phi\|_2^2
\end{equation} 
Therefore we are justified in multiplying the PDEs by $U^{i+1}$ and $w^{i+1}$ respectively and integrating in space.  Indeed, one just chooses a sequence of test functions approximating either $U^{i+1}$ or $w^{i+1}$ and passes the limit through the weak formulation (\ref{stationary:abound}) or (\ref{difference:aweakformulation}), this will be a common technique in the remainder of the work.  In both cases a stronger existence theorem is true but outside the scope of this paper.
\end{Remark}

\bigskip

\subsection{Decay of $w$}

This subsection contains energy decay calculations for $w$, the estimates are an application of the Fourier Splitting Method with bootstrapping.  The first step in the procedure is to apply the Fourier splitting method using the bound (\ref{wenergy:apbound}) to find a preliminary decay rate.  Once established, this preliminary rate is used to deduce a faster decay rate.  This procedure is repeated until the recursion does not lower the rate again, in this case the hold up will be from estimates on the nonlinear term.  The sequence of lemmas leading to Theorem \ref{wildecay:thrm} set up the bootstrap situation which is the main part of the proof for the theorem, establishing (\ref{wi1decay:bound}) is the main goal of this subsection.  The calculations are formal but can be made rigorous by applying them to a sequence of approximating solutions to (\ref{difference:PDE}) (see \cite{MR775190}) or working directly with the weak formulation (see Remark \ref{justifymultiply:remark}).  We begin with an estimate for $|\hat{w}|$.  

\begin{Lemma}
 Let $w^{i+1}$ be the solution of (\ref{difference:aPDE}) given by Theorem \ref{difference:aexist} with $U^i$ and $f$ satisfying the assumptions of the theorem.  Then,
\begin{align}\label{ftw:bound}
 |\hat{w}^{i+1}|\leq C|\xi|\|U^{i}\|_2\left(\int_0^t\|w^{i+1}(s)\|_2\,ds +\nu^{-1}\rho_0^{-1}\|f\|_X\right)
\end{align}
\end{Lemma}
\begin{proof}
 Through the Fourier transform of (\ref{difference:aPDE}), noting the initial data is zero, we write
\begin{equation}
 \hat{w}^{i+1}= -\int_0^t e^{-\nu|\xi|^2(t-s)}(\xi\cdot \widehat{U^{i}w^{i+1}} + \xi\hat{p} + \xi\cdot\widehat{U^{i}\Phi})(s)\, ds\notag
\end{equation}
Young's inequality with the Plancherel theorem bounds
\begin{equation}
 |\widehat{U^{i}w^{i+1}}|+|\widehat{U^{i}\Phi}|\leq \|U^{i}\|_2(\|w^{i+1}\|_2+\|\Phi\|_2)\notag
\end{equation}
Taking the divergence of (\ref{difference:PDE}), then the Fourier transform, one can bound $|\hat{p}|\leq C(|\widehat{U^{i}w^{i+1}}|+|\widehat{U^{i}\Phi}|)$.  All together,
\begin{align}
 |\hat{w}^{i+1}|&\leq C|\xi|\|U^i\|_2\int_0^t(\|w^{i+1}\|_2+\|\Phi\|_2)(s)\, ds\notag\\
&\leq C|\xi|\|U^{i}\|_2\left(\int_0^t\|w^{i+1}(s)\|_2\,ds +\rho_0^{-1}\nu^{-1}\|f\|_X\right)\notag
\end{align}
The last line uses (\ref{assumption1:bound}) to evaluate the integral in time.  This completes the proof.
\end{proof}

\begin{Lemma}\label{bootstrapbase:lemma}
 Let $w^{i+1}$ be the solution of (\ref{difference:aPDE}) given by Theorem \ref{difference:aexist} with $U^i$ and $f$ satisfying the assumptions of the theorem.  Then, for any $m\geq4$, $w^{i+1}$ satisfies the differential inequality
\begin{align}\label{bootstrapbase:ineq}
\frac{d}{dt}&\left((1+t)^m\|w^{i+1}\|_2^2\right)\notag\\
&\qquad \qquad \leq C(m,\rho_0,\nu)\|U^{i}\|^2_2(1+t)^{m-\frac{7}{2}}\left(\int_0^t\|w^{i+1}(s)\|_2\,ds +\|f\|_X\right)^2 \\
&\qquad \qquad \qquad + C\nu^{-3}\|f\|_X^2\|\Phi\|_3^2(1+t)^{m}\notag
\end{align}
\end{Lemma}
\begin{Remark}\label{constant1:rem}
In the statement of the lemma, the constant $C(m,\rho_0,\mu)$ tends to $\infty$ as $\rho_0\rightarrow 0$ or $\nu\rightarrow 0$ and tends toward $0$ as $\nu\rightarrow \infty$.
\end{Remark}
\begin{proof}
Multiply (\ref{difference:aPDE}) by $w^{i+1}$, after integration by parts then application of the bilinear relation (\ref{bilinear:relation}) and the assumed bound (\ref{difference:aass}) we write
\begin{align}\label{step1:difflemma}
\frac{1}{2}\frac{d}{dt}\|w^{i+1}\|_2^2 +\nu\|\nabla w^{i+1}\|_2^2 &= <U^i\cdot\nabla w^{i+1},\Phi>\\
&\leq \|U^i\|_6\|\nabla w^{i+1}\|_2\|\phi\|_3\notag\\
&\leq \frac{C}{\nu^3}\|f\|_X^2\|\Phi\|_3^2 +\frac{\nu}{2}\|\nabla w^{i+1}\|_2^2\notag
\end{align}
Now we \textit{split} the viscous term in Fourier Space around the ball $B(R)$ using the Plancherel theorem:
\begin{align}
-\nu\|\nabla w^{i+1}\|_2^2&\leq -\nu \int_{B(R)^C}|\xi|^2 |\hat{w}^{i+1}|^2\, d\xi \notag\\
&\leq -\nu R^2\int_{B(R)^C} |\hat{w}^{i+1}|^2\, d\xi\notag\\
&\leq-\nu R^2\|\hat{w}^{i+1}\|_2^2 + \nu R^2\int_{B(R)} |\hat{w}^{i+1}|^2\, d\xi \notag
\end{align}
Combining this with  (\ref{step1:difflemma}):
\begin{align}
 \frac{d}{dt}\|w^{i+1}\|_2^2+\nu R^2\|w^{i+1}\|_2^2 &\leq \nu R^2\int_{B(R)}|\hat{w}^{i+1}|\, d\xi + C\nu^{-3}\|f\|_X^2\|\Phi\|_3^2\notag
\end{align}
Then using (\ref{ftw:bound}) we bound
\begin{align}
\int_{B(R)}|&\hat{w}^{i+1}|^2\,d\xi\notag\\
&\leq C\|U^i\|_2^2\left(\int_0^t\|w^{i+1}(s)\|_2\,ds +\nu^{-1}\rho_0^{-1}\|f\|_X\right)\left(\int_{B(R)}|\xi|^2\,d\xi\right)\notag\\
&\leq C\|U^i\|_2^2\left(\int_0^t\|w^{i+1}(s)\|_2\,ds +\nu^{-1}\rho_0^{-1}\|f\|_X\right)R^5\notag
\end{align}
So,
\begin{align}
 \frac{d}{dt}\|w^{i+1}\|_2^2+\nu R^2&\|w^{i+1}\|_2^2 \notag\\
&\leq C\nu R^7 \|U^i\|^2_2(1+\nu^{-1}\rho_0^{-1})^2\left(\int_0^t\|w^{i+1}(s)\|_2\,ds+\|f\|_X\right)^2 \notag\\
& \qquad \qquad + C\nu^{-3}\|f\|_X^2\|\Phi\|_3^2\notag
\end{align}
In the preceding inequality we choose $R^2=\frac{m}{\nu}(1+t)^{-1}$ then use $(1+t)^{m}$ as an integrating factor to establish the lemma.  Examining the line above one can see the constant in the statement of the lemma behaves like $(1+\rho_0^{-1}\nu^{-1})^2\nu^{-7/2}$, this is Remark \ref{constant1:rem}.
\end{proof}

\begin{Theorem}\label{wildecay:thrm}
Let $w^{i+1}$ be the solution of (\ref{difference:PDE}) given by Theorem \ref{difference:aexist} with $U^i$ and $f$ satisfying the assumptions of the theorem.  Then, $w^{i+1}$ satisfies the decay bound
\begin{equation}\label{wi1decay:bound}
\|w^{i+1}(T)\|_2^2 \leq C(\rho_0,\nu)(1+\|U^{i}\|_2^2)^6(1+\|f\|_X^2)\|f\|_X^2(1+T)^{-\frac{5}{2}}
\end{equation}
\end{Theorem}
\begin{Remark}The exponent of $(1+T)$ is such that $\|w^{i+1}\|_2$ is integrable over all time.  The constant $C(\rho_0,\nu)$ in the statement tends to $\infty$ as $\rho_0\rightarrow 0$ or $\nu\rightarrow 0$.  It tends to $0$ as $\nu\rightarrow \infty$ (see Remark \ref{constant1:rem}).
\end{Remark}
\begin{proof}
Combining the bound on $\|w^{i+1}\|_2$ given by (\ref{wenergy:apbound}) with (\ref{bootstrapbase:ineq}) we write
\begin{align}
\frac{d}{dt}\left((1+t)^m\|w^{i+1}\|_2^2\right)&\leq C(m,\rho_0,\nu)\|U^{i}\|^2_2(1+t)^{m-\frac{7}{2}}\left(t^2\|f\|_X^4+\|f\|_X^2\right)\notag\\
&\qquad \qquad +C(\nu)\|f\|_X^2\|\Phi\|_3^2(1+t)^{m}\notag\\
&\leq C(m,\rho_0,\nu)(1+\|U^{i}\|_2^2)(1+\|f\|_X^2)\|f\|_X^2(1+t)^{m-\frac{3}{2}}\notag\\
&\qquad \qquad +C(\nu)\|f\|_X^2\|\Phi\|_3^2(1+t)^{m}\notag
\end{align}
The next step is to integrate in time, the first term on the RHS can be integrated directly while the second term is estimated similar to (\ref{assumption1:heatl3}):
\begin{align}
C(\nu)\|f\|_X^2\int_0^T\|&\Phi(t)\|_3^2(1+t)^{m}\,dt\notag\\ 
&\leq C(\nu)\|f\|_X^2\int_0^T\|\Phi(t)\|_2\|\nabla\Phi(t)\|_2(1+t)^{m}\,dt\notag\\
&\leq C(\nu)\|f\|_X^2\left(\int_0^T (1+t)^m\|\Phi(t)\|_2^2\,ds\right)^{\frac{1}{2}}\left(\int_0^T \|\nabla \Phi(t)\|_2^2\,dt\right)^{\frac{1}{2}}\notag\\
&\leq C(\rho_0,\nu)\|f\|_X^4\notag
\end{align}
This gives an initial decay bound
\begin{equation}\label{bootstrap1:bound}
\|w^{i+1}(T)\|_2^2 \leq C(m,\rho_0,\nu)(1+\|U^{i}\|_2^2)(1+\|f\|_X^2)\|f\|_X^2(1+T)^{-\frac{1}{2}}
\end{equation}

Now we begin the bootstrapping procedure.  Proceeding in a nearly identical way to the argument immediately above, use (\ref{bootstrapbase:ineq}) with (\ref{bootstrap1:bound}) instead of (\ref{wenergy:apbound}), then integrate in time:
\begin{equation}
\|w^{i+1}(T)\|_2^2 \leq C(m,\rho_0,\nu)(1+\|U^{i}\|_2^2)^2(1+\|f\|_X^2)\|f\|_X^2(1+T)^{-1}\notag
\end{equation}
This process can be repeated indefinitely but the ``best'' decay rate will be obtained after six iterations; here ``best'' is meant in the sense of best decay rate obtainable from (\ref{bootstrapbase:ineq}).  That this is in fact the best decay rate can be seen by examining the term $\int_0^t\|w^{i+1}(s)\|_2\,ds$ in (\ref{bootstrapbase:ineq}), once we have established $\|w^{i+1}\|_2 \leq C(1+t)^{-\mu}$ for $\mu>1$ this term integrates to a constant and we obtain the ``best decay rate.''  As the bootstrapping steps are nearly identical to the above arguments and tedious to write out we skip to the final step:
\begin{equation}
\|w^{i+1}(T)\|_2^2 \leq C(m,\rho_0,\nu)(1+\|U^{i}\|_2^2)^6(1+\|f\|_X^2)\|f\|_X^2(1+T)^{-\frac{5}{2}}\notag
\end{equation}
\end{proof}
\begin{Remark}
In the above proof we relies on the exponential decay of $\Phi$ which follows from Assumption \ref{assumption1:ass}.  This can be relaxed to $\hat{f} \sim |\xi|^k$ near the origin for a large $k$.
\end{Remark}

\bigskip

\subsection{Relation between $U^i$ and $w^i$}

In this subsection we make precise, for our approximate solutions, the formal notion $U^i=\int_0^\infty v^i(t)\,dt$.  We show approximations of the integral $\int_0^n v^{i}(t)\,dt$ are bounded uniformly in $L^2$ and are Cauchy with a limit which is a solution of (\ref{stationary:aexist}).  Once this is established we apply the decay results from the previous subsection to find a uniform bound in $L^{2}$ for $U^{i}$.  Throughout this subsection we use $\Phi=e^{\triangle t}f$.

\begin{Lemma}\label{Vcauchy:lemma}
Let $w^{i+1}$ be the solution of (\ref{difference:PDE}) given by Theorem \ref{difference:aexist} with $U^i$ and $f$ satisfying the assumptions of the theorem.  The function $v^{i+1}=w^{i+1}+\Phi$ satisfies $\int_0^\infty v^{i+1}(t)\,dt\in L^2_\sigma$.
\end{Lemma}
\begin{proof}
For each fixed $i$ define the sequence $\{V^{i+1}_n\}_{n\in\mathbb{N}}\subset L^2_\sigma$ by 
\begin{equation}\notag
V^{i+1}_n = \int_0^n v^{i+1}(t)\ dt 
\end{equation}
Since $v^{i+1}(t)\in L^2_\sigma \ a.e.$ the sequence $\{V^{i+1}_n\}$ is well defined.  Relying on Minkowski's inequality for integrals with Assumption \ref{assumption1:ass} (through (\ref{assumption1:bound})) and (\ref{wi1decay:bound}) the following bound shows how the sequence $\{V^{i+1}_n\}$ is bounded uniformly (for $n$) in $L^2_\sigma$:
\begin{align}
\|V_n^{i+1}\|_2&\leq \int_0^n\|v^{i+1}(t)\|_2\,dt\notag\\
&\leq \int_0^n\|w^{i+1}(t)\|_2\,dt +\int_0^n\|\Phi(t)\|_2\,dt\notag\\
&\leq C(\rho_0,\nu)(1+\|U^{i}\|_2^2)^3(1+\|f\|_X^2)^{\frac{1}{2}}\|f\|_X\notag
\end{align}
Similarly, 
\begin{equation}
\|V_{n+1}^{i+1}-V_{n}^{i+1}\|_2 \leq \int_n^{n+1}\|v^{i+1}\|_2\,dt\label{cauchy1}
\end{equation}  
Observing (\ref{wi1decay:bound}) and the decay of $\Phi$ implied by Assumption \ref{assumption1:ass} we know the integral $\int_0^\infty \|v^{i+1}\|_2\, dt$ is finite so the RHS of (\ref{cauchy1}) tends to zero as $n\rightarrow \infty$.  Following well known arguments to prove a contraction lemma we can quickly deduce $\{V_n^{i+1}\}_{n\in\mathbb{N}}$ is Cauchy in $L^2_\sigma$ and has a limit which we label $\int_0^\infty v^{i+1}(t)\,dt$.
\end{proof}
\begin{Remark}
The above lemma also implies $\int_0^\infty v^{i+1}(t)\,dt$ is finite $a.e.$ in $\mathbb{R}^3$.
\end{Remark}

\begin{Lemma}\label{intmakessense:lemma}
Let $w^{i+1}$ be the solution of (\ref{difference:PDE}) given by Theorem \ref{difference:aexist} with $U^i$ and $f$ satisfying the assumptions of the theorem.  The function $v^{i+1}=w^{i+1}+\Phi$ satisfies $\int_0^\infty v^{i+1}(t)\,dt=U^{i+1}$ .
\end{Lemma}
\begin{proof}
To prove this lemma we show $\int_0^\infty v^{i+1}(t)\,dt$ is a weak solution for (\ref{stationary:aPDE}) then use the uniqueness implied by Theorem \ref{stationary:aexist} to conclude the desired result.  Let $\{V_n^{i+1}\}_{n\in\mathbb{N}}$ be as in the previous proof.
%

In (\ref{difference:aweakformulation}) choose $\phi$ to be any member of $\mathcal{V}$ (so that it is constant in time).  Use the relation $v^{i+1}=w^{i+1}+\Phi$ then integrate in time:
\begin{align}
\int_0^n&\left(\frac{d}{dt}<v^{i+1}(t),\phi>+<U^{i}\cdot\nabla v^{i+1}(t),\phi>\right)\,dt\notag\\
 &\qquad\qquad\qquad\qquad\qquad\qquad\qquad=-\nu\int_0^n<\nabla v^{i+1}(t),\nabla\phi>\,dt
\end{align}
After changing the order of integration and evaluating the first integral this becomes
\begin{align}
<v^{i+1}(n),\phi>+<U^{i}\cdot\nabla V_n^{i+1},\phi> &=-\nu<\nabla V_n^{i+1},\nabla\phi>+<f,\phi>\notag
\end{align}
Observe the first term on the LHS tends to zero as $n\rightarrow \infty$.  This follows from the decay bound (\ref{wi1decay:bound}) which implies $w^{i+1}$ tends to zero on compact sets and a similar well known property for the heat equation.  The strong convergence of $\{V_n^{i+1}\}$ in $L^2$ is enough to pass the limit through the remaining terms.  Indeed, if $V^{n+1}=\int_0^\infty v^{i+1}(t)\,dt$ is this limit,
\begin{align}
 |<U^{i}\cdot\nabla (V_n^{i+1}-V^{i+1}),\phi>|\leq C\|\nabla U^{i}\|_2\|V_n^{i+1}-V^{i+1}\|_2\|\nabla \phi\|_3\notag
\end{align}
As $n\rightarrow \infty$ this tends to zero for each test function $\phi\in \mathcal{V}$, hence $V^{i+1}$ is a weak solution of (\ref{stationary:aPDE}).  The uniqueness implied by Theorem \ref{stationary:aexist} finishes the proof of the lemma.
\end{proof}

\begin{Lemma}
Let $U^{i+1}$ be the solution of (\ref{stationary:PDE}) given by Theorem \ref{stationary:aexist} with $U^i$ and $f$ satisfying the assumptions of the theorem.  Then the function $U^{i+1}$ satisfies 
\begin{align}\label{goodl2Ui1:bound}
 \|U^{i+1}\|^2_2\leq C(\rho_0,\nu)(1+\|U^{i}\|_2^2)^6(1+\|f\|_X^2)\|f\|_X^2
\end{align}
\end{Lemma}
\begin{Remark}The constant $C(\rho_0,\nu)$ in the above theorem tends to $\infty$ as $\rho_0\rightarrow 0$ or $\nu\rightarrow 0$.  It tends to $0$ as $\nu\rightarrow \infty$ (see proof of Lemma \ref{bootstrapbase:lemma}).
\end{Remark}
\begin{proof}
Define $v^{i+1}=w^{i+1}+\Phi$.  Just as in the proof of Lemma \ref{Vcauchy:lemma} combine Minkowski's inequality for integral with (\ref{assumption1:bound}) and (\ref{wi1decay:bound}), but this time use the relation from Lemma \ref{intmakessense:lemma}.
\begin{align}
\|U^{i+1}\|_2&\leq \int_0^n\|v^{i+1}(t)\|_2\,dt\notag\\
&\leq \int_0^n\|w^{i+1}(t)\|_2\,dt +\int_0^n\|\Phi(t)\|_2\,dt\notag\\
&\leq C(\rho_0,\nu)(1+\|U^{i}\|_2^2)^3(1+\|f\|_X^2)^{\frac{1}{2}}\|f\|_X\notag
\end{align}
\end{proof}

\bigskip

\subsection{Convergence of $U^{i}$}

The goal of this subsection is to find the limit of the approximating sequence $U^i$ and show this is a solution of the steady state Navier-Stokes equation.  Later in the subsection we make two assumptions on $f$, they are \textit{smallness} assumptions and allow a contraction argument to show $U^i$ is Cauchy.  The assumptions will depend on how \textit{big} we will allow the $L^2_\sigma$ norm of $U$ and throttle $\|\nabla U\|_2$ so that a product of the $L^2_\sigma$ and $\mathring{H}^1_\sigma$ norms of $U$ is small.  We will label this maximum value of the $L^2$ norm $M$ (our choice) and keep it fixed throughout the remainder of this section.

%
\begin{Lemma}\label{ul2M:abound}
Let $U^{i+1}$ be the solution of (\ref{stationary:PDE}) given by Theorem \ref{stationary:aexist} with $U^i$ and $f$ satisfying the assumptions of the theorem.  There exists a constant $C(\rho_0,\nu,M)$ so that if $\|f\|_X\leq C(\rho_0,\nu,M)$ and $\|U^{i}\|_2\leq M$ then $\|U^{i+1}\|_2\leq M$.
\end{Lemma}
\begin{Remark}\label{ul2M:remark}
For fixed $\rho_0$ and $\nu$ the constant $C(\rho_0,\nu,M)$ in the above lemma tends to $0$ as $M\rightarrow \infty$ or $M\rightarrow 0$.  For fixed $M$ the constant tends to $0$ as $\rho_0\rightarrow 0$ or $\nu \rightarrow 0$ and tends to $\infty$ as $\nu \rightarrow \infty$.   
\end{Remark}
\begin{proof}
By setting the RHS of (\ref{goodl2Ui1:bound}) equal to $M^2$ and considering $Z=\|f\|_X^2$ as a variable the proof is reduced to finding roots of the polynomial
\begin{align}
Z^2+Z=L=\frac{M^2}{C(\rho_0,\nu)(1+M^2)^{6}}\notag
\end{align}
Here, $C(\rho_0,\nu)$ is exactly as in (\ref{goodl2Ui1:bound}).  Since $L>0$ This polynomial always has a \textit{strictly positive} root, in this case the root is exactly the constant in the statement of the lemma.  Indeed,
\begin{align}
\frac{-1+\sqrt{1+4L}}{2} \leq \sqrt{L} = \frac{M}{\sqrt{C(\rho_0,\nu)}(1+M^2)^{3}}\label{temp2}
\end{align}
Remark \ref{ul2M:remark} follows by examining the RHS of (\ref{temp2}).
\end{proof}


\begin{Theorem}\label{L2norm1:theorem}
Let $M>0$ and $f$ satisfy Assumption \ref{assumption1:ass}.  There exists a constant $C(\rho_0,\nu,M)$ such that if $\|f\|_X\leq C(\rho_0,\nu,M)$ the following hold:
\begin{itemize}
\item[(i)] The PDE (\ref{stationary:PDE}) has a weak solution $U\in H^1_\sigma$ (in the sense of (\ref{weakstatement:stationary})).
\item[(ii)] This solution satisfies $\|U\|_2\leq M$ in addition to (\ref{stationary:bound})
\item[(iii)] This solution is unique among all solutions which satisfy (\ref{stationary:bound}) and have a finite $L^2_\sigma$ norm.
\end{itemize}
\end{Theorem}
\begin{Remark}
For fixed $\rho_0$ and $\nu$ the constant $C(\rho_0,\nu,M)$ in the above lemma tends to $0$ as $M\rightarrow \infty$ or $M\rightarrow 0$.  For fixed $M$ the constant tends to $0$ as $\rho_0\rightarrow 0$ or $\nu \rightarrow 0$ and tends to $\infty$ as $\nu \rightarrow \infty$.  The behavior of the constant as $M\rightarrow \infty$ with the bound (\ref{stationary:bound}) implies $\|U\|_2\|\nabla U\|_2 \leq C(\nu)$.  In the time dependent Navier-Stokes system it is well known that when this product of norms is small for initial data the solution will remain smooth and well behaved, this result fits into that regime.
\end{Remark}
\begin{proof}
Choose $U^0\in H^1_\sigma$ so that $\|U^0\|_2\leq M$ and $\|\nabla U^0\|_2^2\leq \nu^{-2}\|f\|_X^2$.  To construct such a function one could fix $f$ then take a solution $\bar{U}$ for (\ref{stationary:PDE}).  At this point the solution is not known to be unique or have finite $L^2_\sigma$ norm but by using a suitable cut-off function in Fourier space ($U^0=\check{\chi}\ast \bar{U}$ where $\chi$ is equal to zero inside a ball containing the origin and one elsewhere) it is possible to limit the $L^2_\sigma$ norm while not increasing the $\mathring{H}^1_\sigma$ norm.

Starting with $U^0$, solve (\ref{stationary:aPDE}) recursively using Theorem \ref{stationary:aexist} to find a sequence $\{U^i\}_{i=0}^\infty$ which satisfies $\|\nabla U^i\|_2\leq \nu^{-2}\|f\|_X^2$.  Lemma \ref{ul2M:abound} proviedes the uniform bound $\|U^i\|_2\leq M$ and so its limit, if it exists, must also satisfy this bound.  We will now show this sequence is Cauchy in $\mathring{H}^1_\sigma$ and a limit does indeed exist.  The difference $Y^{i+1}=U^{i+1}-U^i$ solves
\begin{align}\label{stationarydifference:approxPDE}
U^{i}\cdot\nabla Y^{i+1} + Y^i\cdot\nabla U^{i}+\nabla p &=\nu\triangle Y^{i+1}\notag
\end{align}
After multiplying this by $Y^{i+1}$, integrating by parts and using the bilinear relation (\ref{bilinear:relation}) one can deduce
\begin{align}
\nu \|\nabla Y^{i+1}\|_2^2 &= <Y^i\cdot\nabla Y^{i+1},U>\notag\\
&\leq \|Y^i\|_6\|\nabla Y^{i+1}\|_2\|U\|_3\notag\\
&\leq \frac{1}{2\nu}\|\nabla Y^i\|_2^2\|U\|_2\|\nabla U\|_2 + \frac{\nu}{2}\|\nabla Y^{i+1}\|_2^2\notag
\end{align}
The above sequence relies on H\"older's inequality, Cauchy's inequality, and the Gagliardo-Nirenberg-Sobolev inequality.  It implies
\begin{align}
 \|\nabla Y^{i+1}\|_2^2 &\leq C\nu^{-2}\|U^i\|_2\|\nabla U^i\|_2\|\nabla Y^i\|_2^2\notag\\
&\leq C\nu^{-3}M\|f\|_X\|\nabla Y^i\|_2^2
\end{align}
Note the multiplication by $Y^{i+1}$ is justified since all $U^{i}$ (and hence all $Y^{i}$) are bounded in $H^1_\sigma$.
Using this bound recursively one finds
\begin{align}
 \|\nabla Y^{i+1}\|_2^2 &\leq (C\nu^{-3}M\|f\|_X)^{i+1}\|\nabla Y^1\|_2^2\notag\\
&\leq (C\nu^{-3}M\|f\|_X)^{i+1}2\nu^{-2}\|f\|_X^2\notag
\end{align}
The last step relies the uniform bound on $\|\nabla U^i\|_2\leq \nu^{-2}\|f\|_X^2$.  If $\|f\|_X< \frac{\nu^3}{CM}$ where $C$ is the same as the line immediately above then $Y^i$ tends to zero in $\mathring{H}^1_\sigma$ which implies $U^i$ is Cauchy, call its limit $\tilde{U}$.  Through this construction we can also be sure $\|\tilde{U}\|_2\leq M$.  Using standards arguments this strong convergence is enough to pass the limit through (\ref{stationary:aPDE}) and show $\tilde{U}$ is a solution of (\ref{stationary:PDE}).  For completeness we will demonstrate how to pass through the nonlinear term:
\begin{align}
|<\tilde{U}\cdot\nabla \tilde{U}, \phi> - <U^i\cdot\nabla U^{i+1}, \phi>| &\leq I+II\notag
\end{align}
\begin{align}
I&= |<\tilde{U}\cdot\nabla (\tilde{U}-U^{i+1}),\phi>|\notag\\
II&= |<(\tilde{U}-U^{i})\cdot\nabla U^{i+1},\phi>|\notag
\end{align}
To show $I\rightarrow 0$ use H\"older's inequality:
\begin{align}
 I \leq \|\tilde{U}\|_3\|\nabla (\tilde{U}-U^{i+1})\|_2\|\phi\|_6\notag
\end{align}
Since the $L^3$ norm of $\tilde{U}$ and the $L^6$ norm of $\phi$ are bounded, the strong convergence $U^{i}\rightarrow \tilde{U}$ in $\mathring{H}^1_0$ shows the RHS tends to zero.  The term $II$ is handled in a nearly identical way.

It remains to establish that $\tilde{U}$ is the unique solution of (\ref{stationary:PDE}) among all solutions which satisfy (\ref{stationary:bound}) and have finite $L^2$ norm.  Let $U$ be any other solution which satisfies (\ref{stationary:bound}) and has a finite $L^2_\sigma$ norm.  The difference $Y=U-\tilde{U}$ solves 
\begin{align}
U\cdot\nabla Y + Y\cdot\nabla \tilde{U}+\nabla p &=\nu\triangle Y\label{DIFF:eq}
\end{align}
$U$ and $\tilde{U}$ are bounded in $L^2_\sigma$ and $\mathring{H}^1_\sigma$ we are allowed to multiply this equation by $Y$.  Then, proceeding in the same way as the lines leading to (\ref{stationarydifference:approxPDE}),
\begin{align}
 \|\nabla Y\|_2^2 & \leq C\nu^{-3}M\|f\|_X\|\nabla Y\|_2^2\notag
\end{align}
The assumption on $f$ made earlier in this proof is enough to guarantee  $\|f\|_X<\nu^3/CM$ and implies the solution is unique.
\end{proof}
\begin{Remark}
This is exactly Theorem \ref{L2norm:theorem}.  Examining (\ref{DIFF:eq}) it seems $Y\in \mathring{H}^1_\sigma$ might be enough to obtain uniqueness as the first term on the LHS would formally integrate to zero after multiplying by $Y$.  Unfortunately, following techniques used in this paper, we are not able to multiply (\ref{DIFF:eq}) by $Y$ unless we know also $Y\in L^3_\sigma$.  Indeed, using H\"older's inequality and the Gagliardo-Nirenberg-Sobolev inequality one can see $U\cdot\nabla Y\in (\mathring{H}^1_\sigma\cap L^3_\sigma)'\subset(H^1_\sigma)'$, but not $U\cdot\nabla Y\in (\mathring{H}^1_\sigma)'$.

The uniqueness in this theorem, with the exact same proof, could instead be stated ``This solution is unique among all solutions which satisfy (\ref{stationary:bound}) and have a finite $L^3$ norm'' and using some other technique it \textit{may} be possible to expand this uniqueness theorem further.
\end{Remark}

\bigskip

\section{Stability of Solutions}
An important property of steady state solutions for physical problems is stability: ``If a steady state solution is perturbed will it return to the same solution?''  In the setting of the Navier-Stokes equation we investigate the stability of a solution $U$ for (\ref{stationary:PDE}) by considering a perturbation $w_0$ and examining the long term behavior of the solution of the Navier-Stokes equation (\ref{NS:PDE}) with initial data $u(0)= U+w_0$.
In particular, one would like to know what conditions on $w_0$ will guarantee solutions of (\ref{NS:PDE}) approach $U$ as time becomes large.  An equivalent problem, found by subtracting (\ref{stationary:PDE}) from the (\ref{NS:PDE}), is to determine when solutions of the following PDE tend to $0$:
\begin{align}\label{perturbeddifference:PDE}
w_t+u\cdot\nabla w + w\cdot\nabla U+\nabla p &=\nu\triangle w\\
\nabla \cdot w=0 \ \ \ w(0)&=w_0\notag
\end{align}
A through examination of stability for the steady state Navier-Stokes equation is currently outside the reach of modern techniques (even in bounded domains) but using energy techniques we can prove strong stability results for perturbations of finite energy under certain restraints on $f$.  Following the literature we introduce the following notion of stability, commonly called nonlinear stability, and find conditions on $f$ which guarantee this type of stability.

\begin{Definition}\label{NLstability:dfn}
We say a solution $U\in H^1_\sigma$ of (\ref{stationary:PDE}) corresponding to $f\in X$ and given by Theorem \ref{L2norm:theorem} is \textit{nonlinearly stable} if it satisfies the following: if $w_0\in L^2_\sigma$ is a perturbation and $u$ is the solution of the Navier-Stokes equation (\ref{NS:PDE}) with initial data $w_0+U$ which satisfies, for any $T>0$,
\begin{equation}
u\in L^{\infty}(0,T;L^2)\cap L^{2}(0,T;\mathring{H}^1_\sigma \notag
\end{equation}
then
\begin{itemize}
 \item[(i)] For every $\epsilon >0$ there is a $\delta>0$ such that 
\begin{equation}\notag
 \|w_0\|_2\leq \delta\ \ \  \mathrm{implies} \ \ \ \sup_{t\in\mathbb{R}^+}\|u(t)- U\|_2 \leq \epsilon
\end{equation}
 \item[(ii)]$u(t)$ tends to $U$ as time becomes large, that is
\begin{equation}\notag
 \lim_{t\rightarrow\infty}\|u(t)-U\|_2=0
\end{equation}
\end{itemize}
\end{Definition}

To start we can multiply (\ref{perturbeddifference:PDE}) by $w$, then integrate by parts and use the bilinear relation (\ref{bilinear:relation}) to find a formal energy inequality
\begin{align}
\frac{1}{2}\frac{d}{dt}\|w\|_2^2+\nu\|\nabla w\|_2^2 &=-<w\cdot\nabla U, w>\notag\\
&\leq \|w\|_6\|U\|_3\|\nabla w\|_2\notag\\
&\leq C\|U\|^{\frac{1}{2}}_2\|\nabla U\|^{\frac{1}{2}}_2\|\nabla w\|^2_2\notag
\end{align}
The last two lines above were obtained with a combination of H\"older's inequality, Cauchy's inequality, and the Gagliardo-Nirenberg-Sobolev inequality.
If $f$ is chosen so that $U$ is small in $H^1_\sigma$ norm ($C\|U\|_2^{1/2}\|\nabla U\|_2^{1/2}\leq \nu/2$) the above inequality becomes
\begin{align}\label{formalenergy:ineq}
\frac{d}{dt}\|w\|_2^2+\nu\|\nabla w\|_2^2 &\leq 0
\end{align}
This differential inequality implies $\|w\|_2^2$ is bounded uniformly.  As $U\in L^2_\sigma$ we can say the same about $u=w+U$.  In other words, if one considers a solution $U$ given by Theorem \ref{L2norm:theorem} then any finite energy perturbation will stay close to $U$, this is condition (i) in Definition \ref{NLstability:dfn}.  The rest of this section will be spent proving the stronger statement, that all finite energy perturbations return to $U$.

\bigskip

\subsection{Existence Theorems}
Here we state the existence theorems and properties of the two PDEs examined in this section.  Proofs of these theorems can be found in literature and are omitted.
\begin{Theorem}\label{perturbed:theorem}
Given $T>0$, initial data $u_0\in L^2_\sigma$, and a forcing function $f\in \mathring{H}^1_\sigma$ the PDE (\ref{NS:PDE}) has a weak solution
\begin{equation}\label{perturbed:inspace}
u\in L^{\infty}([0,T];L^2_\sigma)\cap L^2([0,T];H^1_\sigma)
\end{equation}
in the whole space $\mathbb{R}^3$ which satisfies, for any $\phi \in \mathcal{V}$,
\begin{equation}
 <u_t,\phi> + \nu<\nabla u,\nabla \phi> + <u\cdot\nabla u,\phi> = <f,\phi>\notag
\end{equation}

\end{Theorem}
\begin{proof}
See \cite{MR673830}, \cite{MR972259}, \cite{MR589434}, \cite{MR0050423}, \cite{MR1555394}, and \cite{MR1846644}.
\end{proof}

\begin{Theorem}\label{perturbeddifferenece:theorem}
Let $u$ satisfy (\ref{perturbed:inspace}) and $U\in H^1_\sigma$.  There is a constant $C(\nu)$ such that if $\|U\|_{H^1_\sigma}\leq C(\nu)$ the PDE (\ref{perturbeddifference:PDE}) has a unique weak solution 
\begin{equation}\notag
w\in L^\infty(\mathbb{R}^+,L^2_\sigma)\cap L^2(\mathbb{R}^+,\mathring{H}^1_\sigma)
\end{equation}
satisfying, for any $\phi\in \mathcal{V}$,
\begin{equation}
 <w_t,\phi> + \nu<\nabla w,\nabla \phi> + <u\cdot\nabla w,\phi> + <w\cdot\nabla U,\phi> = <f,\phi>\notag
\end{equation}
 as well as the following energy inequalities:
\begin{align}
\frac{d}{dt}\|w\|_2^2+\nu\|\nabla w\|_2^2 &\leq 0\label{perturbeddifferenece:difineq}\\
\sup_{t\in\mathbb{R}^+} \|w(t)\|_2^2+\nu\int_0^\infty \|\nabla w\|_2^2 &\leq \|w(0)\|_2^2\label{perturbeddifference:ineq}
\end{align}
\end{Theorem}
\begin{proof}
The proof of this theorem is similar to the one immediately preceding with exception of the a priori bounds which are argued formally preceding (\ref{formalenergy:ineq}).  The linearity of the equation implies uniqueness (which is not known for general solutions of (\ref{NS:PDE})).
\end{proof}

\bigskip

\subsection{Decay of $w$}
In this subsection we prove a decay property for solutions $w$ given by Theorem \ref{perturbeddifferenece:theorem}.  We show $\|w\|_2\rightarrow 0$ using a method developed for the Navier-Stokes equation in \cite{MR1432588}.  The method relies on generalized energy inequalities for the solutions which allow the energy to be decomposed in high and low frequencies.  These are estimated independently and shown to approach zero.

\begin{Lemma}
Let $u$ and $U$ satisfy the assumptions of Theorem \ref{perturbeddifferenece:theorem} with $\|U\|_{H^1_\sigma}$ less then the given constant.  Let $\phi=e^{-|\xi|^2}$, $\psi=1-\phi$ and $E(t)\in C^1([0,\infty);L^\infty)$.  The solution given by Theorem \ref{perturbeddifferenece:theorem} satisfies the following two generalized energy inequalities: 
\begin{align}
\|\check{\phi}\ast w(t)\|_2^2&\leq\|e^{\nu\triangle(t-s)}\check{\phi}\ast w(s)\|_2^2\notag\\
&\qquad +2\int_s^t |<u\cdot\nabla w,e^{2\nu\triangle (t-\tau)}\check{\phi}\ast\check{\phi}\ast w>|\, d\tau\notag\\
&\qquad \qquad+2\int_s^t |<w\cdot\nabla U,e^{2\nu\triangle (t-\tau)}\check{\phi}\ast\check{\phi}\ast w>|\, d\tau\label{perturbeddifference:geq1}\\
E(t)\|\psi \hat{w}(t)\|_2^2&\leq E(s)\|\psi \hat{w}(s)\|_2^2 \notag\\
&\qquad-2\nu\int_s^tE(\tau)\|\xi\psi\hat{w}(\tau)\|_2^2\, d\tau
+\int_s^tE'(\tau)\|\psi\hat{w}(\tau)\|_2^2\, d\tau \notag\\ 
&\qquad\qquad+2\int_s^tE(\tau) |<\widehat{w\cdot\nabla U},\psi^2 \hat{w}>|\, d\tau\label{perturbeddifference:geq2}\\
&\qquad\qquad\qquad+2\int_s^tE(\tau) |<\widehat{u\cdot\nabla w},(1-\psi^2)\hat{w}>|\, d\tau\notag
\end{align}
\end{Lemma}
\begin{proof}
We give a formal proof here which can be made precise by considering an approximating sequence, see \cite{MR1432588} for details.

To see the first inequality multiply the PDE (\ref{perturbeddifference:PDE}) by $e^{2\nu\triangle(t+s)}\check{\phi}\ast\check{\phi}\ast w$ and integrate from $s$ to $t$.  The assumptions are enough to ensure all integrals are finite and this multiplication makes sense.  After integration by parts:
\begin{align}
\|\check{\phi}\ast w(t)\|_2^2&\leq\|e^{\nu\triangle(t-s)}\check{\phi}\ast w(s)\|_2^2 - \nu\int_s^t\|\nabla(e^{\nu\triangle(t-\tau)}\check{\phi})\ast w\|_2^2\, d\tau\notag\\
&\qquad+\int_s^t<\partial_{\tau}(e^{\nu\triangle(t-\tau)}\check{\phi})\ast w, e^{\nu\triangle(t-\tau)}\check{\phi}\ast w>\,d\tau\notag\\
&\qquad \qquad+2\int_s^t |<u\cdot\nabla w,e^{2\nu\triangle (t-\tau)}\check{\phi}\ast\check{\phi}\ast w>|\, d\tau\notag\\
&\qquad \qquad\qquad +2\int_s^t |<w\cdot\nabla U,e^{2\nu\triangle (t-\tau)}\check{\phi}\ast\check{\phi}\ast w>|\, d\tau\notag
\end{align}
$e^{\nu\triangle(t+s)}\check{\phi}$ describes a heat flow so the second and third terms on the RHS add to zero, this proves (\ref{perturbeddifference:geq1}).  For the second inequality, take the Fourier Transform of (\ref{perturbeddifference:PDE}) then multiply by $\psi^2\hat{w}$.  After integration by parts one finds 
\begin{align}E(t)\|\psi \hat{w}(t)\|_2^2&\leq E(s)\|\psi \hat{w}(s)\|_2^2 -2\nu\int_s^tE(\tau)\|\xi\psi\hat{w}(\tau)\|_2^2\, d\tau\notag\\
&\qquad+\int_s^tE'(\tau)\|\psi\hat{w}(\tau)\|_2^2\, d\tau +2\int_s^tE(\tau) |<\widehat{w\cdot\nabla U},\psi^2 \hat{w}>|\, d\tau\notag\\
&\qquad\qquad+2\int_s^tE(\tau) |<\widehat{u\cdot\nabla w},\psi^2\hat{w}>|\, d\tau\notag
\end{align}
The bilinear relation (\ref{bilinear:relation}) and the Plancherel theorem imply $<\widehat{u\cdot\nabla w},\hat{w}>=0$ and (\ref{perturbeddifference:geq2}) follows immediately.
\end{proof}

\begin{Theorem}\label{perturbeddifferencedecay:theorem}
Let $u$ and $U$ satisfy the assumptions of Theorem \ref{perturbeddifferenece:theorem} with $\|U\|_{H^1_\sigma}$ less then the given constant and $\sup_t\|u\|_2 <\infty$.  The energy of the solution given by Theorem \ref{perturbeddifferenece:theorem} decays to zero.  That is,
\begin{equation}
\lim_{t\rightarrow 0}\|w(t)\|_2=0
\end{equation}
\end{Theorem}
\begin{proof}
Following \cite{MR1432588} we bound first the low frequencies using (\ref{perturbeddifference:geq1}) and then the high frequencies using (\ref{perturbeddifference:geq2}) and the Fourier Splitting Method.

To show the low frequencies tend to zero we start by estimating the integrals on the RHS of (\ref{perturbeddifference:geq1}):
\begin{align}
|<u\cdot\nabla w,e^{2\nu\triangle (t-\tau)}\check{\phi}\ast\check{\phi}\ast w>| &= | <\check{\phi}\ast\check{\phi}\ast u\cdot w,e^{2\nu\triangle (t-\tau)}\nabla w>|\notag\\
&\leq \|\check{\phi}\ast\check{\phi}\ast u\cdot w\|_2\|\nabla w\|_2\notag\\
&\leq C\|u\|_2\|\nabla w\|^2_2\notag
\end{align}
This estimate was obtained using integration by parts, the Cauchy-Schwartz inequality, Young's inequality, and the Gagliardo-Nirenberg-Sobolev inequality.  Similarly,
\begin{align}
|<w\cdot\nabla U,e^{2\nu\triangle (t-\tau)}\check{\phi}\ast\check{\phi}\ast w>|\leq C\|U\|_2\|\nabla w\|^2_2\notag
\end{align}
Combining these two bounds with (\ref{perturbeddifference:geq1}) yields:
\begin{align}
\|\check{\phi}\ast w(t)\|_2^2&\leq\|e^{\nu\triangle(t-s)}\check{\phi}\ast w(s)\|_2^2 +C\left(\sup_{\tau\in \mathbb{R}^+}\|u(\tau)\|_2^2+\|U\|_2^2\right)\int_s^t \|\nabla w\|_2^2\, d\tau\notag
\end{align}
Heat energy is known to approach zero as time becomes large so
\begin{align}
\limsup_{t\rightarrow\infty}\|\check{\phi}\ast w(t)\|_2^2&\leq C\left(\sup_{\tau\in \mathbb{R}^+}\|u(\tau)\|_2^2+\|U\|_2^2\right)\int_s^\infty \|\nabla w\|_2^2\, d\tau\notag
\end{align}
The LHS is independent of $s$, noting the energy bound (\ref{perturbeddifference:ineq}) we see the RHS tends to zero as $s\rightarrow \infty$.  Using the Plancherel theorem we conclude
\begin{align}\label{lowfreqbound}
\lim_{t\rightarrow\infty}\|\phi\hat{w}(t)\|_2^2=\lim_{t\rightarrow\infty}\|\check{\phi}\ast w(t)\|_2^2=0
\end{align}

We begin work with the high frequencies on a similar path, bounding the integrals on the RHS of (\ref{perturbeddifference:geq2}).  Note $\psi=1-e^{-|\xi|^2}\in L^\infty$, then
\begin{align}
|<\widehat{w\cdot\nabla U},\psi^2 \hat{w}>|&=|<\xi \cdot \widehat{w\cdot U},\psi^2 \hat{w}>|\notag\\
&\leq \|\widehat{w\cdot U}\|_2\|\xi \hat{w}\|_2\notag\\
&\leq C\|w\|_6\|U\|_3\|\nabla w\|_2\notag\\
&\leq C\|U\|_3\|\nabla w\|^2_2\notag
\end{align}
This chain of inequalities used the Cauchy-Schwartz inequality, the Plancherel theorem, H\"older's inequality, and the Gagliardo-Nirenberg-Sobolev inequality.  Similarly, but this time making use of the rapid decay properties of $1-\psi^2$,
\begin{align}
|<\widehat{u\cdot\nabla w},(1-\psi^2)\hat{w}>|&\leq \|(1-\psi^2)\widehat{u\cdot w}\|_2\|\xi\hat{w}\|_2\notag\\
&\leq C\|(1-\psi^2)^\vee\|_{6/5}\|u\|_2\|w\|_6\|\nabla w\|_2\notag\\
&\leq C\|u\|_2\|\nabla w\|^2_2\notag
\end{align}
Use these two bounds with (\ref{perturbeddifference:geq2}) to find 
\begin{align}\label{perturbeddifference:geq22}
\|\psi \hat{w}(t)\|_2^2&\leq \frac{E(s)}{E(t)}\|\psi \hat{w}(s)\|_2^2\notag\\ &\qquad-2\nu\int_s^t\frac{E(\tau)}{E(t)}\|\xi\psi\hat{w}(\tau)\|_2^2\, d\tau+\int_s^t\frac{E'(\tau)}{E(t)}\|\psi\hat{w}(\tau)\|_2^2\, d\tau\notag\\
&\qquad\qquad +C\left(\sup_{\tau\in\mathbb{R}^+}\|u(\tau)\|_2^2+\|U\|_2^2\right)\int_s^t\frac{E(\tau)}{E(t)} \|\nabla w\|_2^2\, d\tau
\end{align}
Now \textit{split} the viscous term  and the term with $E'$ around the ball with radius $\rho(\tau)>0$, $B(\rho)$:
\begin{align}
-2\nu\int_s^t\frac{E(\tau)}{E(t)}\|\xi\psi&\hat{w}(\tau)\|_2^2\, d\tau+\int_s^t\frac{E'(\tau)}{E(t)}\|\psi\hat{w}(\tau)\|_2^2\, d\tau\notag\\ 
&\leq\frac{1}{E(t)}\int_s^t(E'(\tau)-2\nu E(\tau)\rho(\tau)^2)\int_{B(\rho)^C}|(1-\phi)\hat{w}|^2\, d\xi\, d\tau\notag\\
&\qquad +\int_s^t\frac{E'(\tau)}{E(t)}\int_{B(\rho)}|(1-\phi)\hat{w}|^2\,d\xi\, d\tau\notag
\end{align}
Upon choosing $E(\tau)=(1+t)^\alpha$ and $\rho^2=\alpha/2\nu(1+t)$ ($\alpha >3$), (\ref{perturbeddifference:geq22}) becomes
\begin{align}
\|(1-\phi) \hat{w}(t)\|_2^2&\leq \frac{(1+s)^\alpha}{(1+t)^\alpha}\|(1-\phi) \hat{w}(s)\|_2^2\notag\\ &\qquad +\int_s^t\frac{\alpha(1+\tau)^{\alpha-1}}{(1+t)^\alpha}\int_{B(\rho)}|(1-\phi)\hat{w}|^2\,d\xi\, d\tau\notag\\
&\qquad \qquad +C\left(\sup_{\tau\in\mathbb{R}^+}\|u(\tau)\|_2^2+\|U\|_2^2\right)\int_s^t \|\nabla w\|_2^2\, d\tau\notag
\end{align}
Note $|1-\phi|\leq |\xi|^2$ if $|\xi|<1$, so for large values of $s$,
\begin{align}
\int_{B(\rho)}|(1-\phi)\hat{w}|^2\,d\xi\leq C(1+\tau)^{-2}\|w(\tau)\|_2^2\notag
\end{align}
This implies, again for large $s$,
\begin{align}
\int_s^t\frac{\alpha(1+\tau)^{\alpha-1}}{(1+t)^\alpha}\int_{B(\rho)}|(1-\phi)\hat{w}|^2\,d\xi\, d\tau
&\leq C\sup_{\tau\in \mathbb{R}^+}\|w(\tau)\|_2^2\int_s^t\frac{(1+\tau)^{\alpha-3}}{(1+t)^\alpha}\, d\tau\notag\\
&\leq C\sup_{\tau \in \mathbb{R}^+}\|w(\tau)\|_2^2(1+t)^{-2}
\end{align}
Taking into account the energy bound (\ref{perturbeddifference:ineq}), this tends to zero as $t$ becomes large.  For any large $s$ we are now justified in writing
\begin{align}
\limsup_{t\rightarrow\infty}\|(1-\phi) \hat{w}(t)\|_2^2&\leq C\left(\sup_{\tau\in\mathbb{R}^+}\|u(\tau)\|_2^2+\|U\|_2^2\right)\int_s^\infty \|\nabla w\|_2^2\, d\tau\notag
\end{align}
Again relying on (\ref{perturbeddifference:ineq}) then letting $s\rightarrow \infty$ we find
\begin{align}
\lim_{t\rightarrow\infty}\|(1-\phi) \hat{w}(t)\|_2^2=0\notag
\end{align}
Using this limit in the triangle inequality with the low frequency limit (\ref{lowfreqbound}) completes the proof.
\end{proof}

\begin{Theorem}
Let $f$ satisfy the assumptions of Theorem \ref{L2norm:theorem} and be such that $\|f\|_X$ is less then the constant given by the theorem and $\|U\|_{H^1_\sigma}$ is less then the constant given by Theorem \ref{perturbeddifferenece:theorem}.  The solution $U$ of (\ref{stationary:PDE}) is nonlinearly stable in the sense of Definition \ref{NLstability:dfn}.
\end{Theorem}
\begin{proof}
Let $u$ be given by Theorem \ref{perturbed:theorem}, the difference $v=u-U$  solves (\ref{perturbeddifference:PDE}) and $u$, $U$, $f$ meet the criteria of Theorem \ref{perturbeddifferencedecay:theorem}.  This proves (ii) in Definition \ref{NLstability:dfn}.  Integrating (\ref{perturbeddifferenece:difineq}) in time proves (i) in Definition \ref{NLstability:dfn}.
\end{proof}


\end{document}